\documentclass{IEEEtran}
\usepackage{amsmath,amsfonts,amssymb}
\usepackage{multicol,multirow,diagbox,bm}
\usepackage{graphicx,subfigure}
\usepackage{lineno,cases}
\allowdisplaybreaks 
\usepackage{pgf,tikz}
\usetikzlibrary{arrows,shapes,automata,positioning}
\usetikzlibrary{fit}
\usepackage{pst-node}
\usepackage{mathptmx}
\usepackage{cite}
\usepackage{color,xcolor}
\renewcommand{\bar}{\overline}
\renewcommand{\top}{\intercal}
\newcommand{\lm}{\lambda}
\newcommand{\e}{\epsilon}
\newcommand{\R}{\mathbb{R}}

\newcommand{\W}{\mathbb{W}}
\newtheorem{assj}{Assumption}
\newtheorem{lemj}{Lemma}
\newtheorem{thmj}{Theorem}
\newtheorem{remj}{Remark}
\newtheorem{prop}{Problem}
\newtheorem{defn}{Definition}
\newcommand{\pb}{\begin{IEEEproof} }
\newcommand{\pe}{\end{IEEEproof}}

\usepackage{cite}
\begin{document}
	\title{
		Optimal Output Consensus for  Nonlinear Multi-agent Systems with Both Static and Dynamic Uncertainties \thanks{This work was supported in part by National Natural Science Foundation of China under Grants 61973043 and 61873250 and in part by USTC Research Funds of the Double First-Class Initiative under Grant YD2100002002. ({\it Corresponding author: Xinghu Wang)}}
	}
	\author{Yutao Tang and Xinghu Wang
		\thanks{Y. Tang is with the School of Automation, Beijing University of Posts and Telecommunications, Beijing 100876, China (e-mail: yttang@bupt.edu.cn). X. Wang is with Department of Automation, University of Science and Technology of China, Hefei, 230027, China (e-mail: xinghuw@ustc.edu.cn). }
	}
	\maketitle

\begin{abstract}
	In this technical note, we investigate an optimal output consensus problem for heterogeneous uncertain nonlinear multi-agent systems. The considered agents are described by high-order nonlinear dynamics subject to both static and dynamic uncertainties. A two-step design, comprising sequential constructions of optimal signal generator and distributed partial stabilization feedback controller, is developed to overcome the difficulties brought by nonlinearities, uncertainties, and optimal requirements. Our study can not only assure an output consensus, but also achieve an optimal agreement characterized by a distributed optimization problem. 
\end{abstract}
	
	
\IEEEpeerreviewmaketitle
	
\section{Introduction}

In the past few years, distributed optimization has attracted much attention due to its broad potential applications in multi-robot systems, smart grid and sensor networks. In a typical setting,  each agent has access to a private objective function and all agents are regulated to achieve a consensus on the optimal solution of the sum of all local functions. Many important results were obtained based on gradients or subgradients of the local objective functions combined with consensus rules, including both discrete-time and continuous-time algorithms  \cite{nedic2010constrained,shi2013reaching,jakovetic2015linear,yang2017distributed,zeng2017distributed,li2018distributed}. 

Since distributed optimization tasks may be implemented or depend on  physical dynamics in practice,  optimal consensus involving high-order agent dynamics deserves further investigation. Compared with the pure (output) consensus problem, the consensus point for all outputs of agents is additionally required to be an optimal solution of the global cost function. Note that this optimal solution can only be determined and reached in a distributed way.  Some interesting attempts have been made in  \cite{zhang2017distributed,xie2019global,qiu2019distributed} for integrator agents, \cite{tang2017cyb} for linear agents, and \cite{wang2016distributed,tang2018ijrnc} for special classes of nonlinear agents. However, optimal output consensus for more general nonlinear multi-agent systems is still far from being solved, especially for agents being heterogeneous and subject to uncertainties.  

In this paper, we consider nonlinear multi-agent systems in the Byrnes-Isidori normal form which can model many typical mechanical and electromechanical systems \cite{khalil2002nonlinear}.  In literature, there have been many consensus results for agents of this type, e.g., \cite{kim2011output,liu2013distributed, rezaee2015average}. This normal form is general enough to cover the dynamics reported in existing optimal consensus results\cite{kia2015distributed, gharesifard2014distributed, zhang2017distributed, tang2017cyb, wang2016distributed, tang2018ijrnc,xie2019global,qiu2019distributed}. Here, we further take into account heterogeneous nonlinear dynamics having both static and dynamic uncertainties, which inevitably bring technical difficulties in resolving the optimal output consensus problem. In a preliminary work \cite{tang2018acc}, this problem was studied for such class of agents assuming that the compact set containing static uncertainties is prior known. In this present study, we remove such restrictive condition and allow the boundary of this compact set to be unknown.  

The contribution of this paper is at least two-fold.  First, we solve the optimal output consensus problem for a larger class of uncertain nonlinear multi-agent systems, significantly improving the existing results reported in \cite{tang2017cyb,wang2016distributed,tang2018ijrnc,qiu2019distributed}.  Second, a novel dynamic compensator based distributed controller is developed for  effectively addressing complicated uncertainties, while precise information of system dynamics is required in \cite{zhang2017distributed, xie2019global, qiu2019distributed}. Moreover, in contrast with relevant results in \cite{wang2016distributed,tang2018acc}, the boundary of the compact set containing uncertain parameters is allowed to be unknown. 

The rest of this paper is organized as follows. Preliminaries and problem formulation are presented in Section \ref{sec:formulation}. Then, the design scheme and main results are provided in Sections \ref{sec:two} and \ref{sec:main} with an illustrative example in Section \ref{sec:simu}. Finally, conclusions are given in Section \ref{sec:con}.
 
{\em Notation:} Let $\R^N$ be the $N$-dimensional Euclidean space. Denote $\mbox{col}(a_1,\, \dots,\,a_N) ={[a_1^\top,\, \dots,\,a_N^\top]}^\top$ for vectors $a_1,\,\dots,\,a_N$. ${\bm 1}_N$ (or ${\bm 0}_N$) denotes an all-one (or all-zero) vector in $\R^N$ and ${I}_N$ denotes the $N \times N$ identity matrix. Let $M_1=\frac{1}{\sqrt{N}}{\bm 1}_N$ and $M_2$ be the matrices satisfying $M_2^\top M_1={\bm 0}_{N-1}$, $M_2^\top M_2=I_{N-1}$, and $M_2 M_2^\top=I_{N}-M_1 M_1^\top$. Denote the Euclidean norm of vector $a$ by $||a||$ and the spectral norm of matrix $A$ by $||A||$. A continuous function $\alpha\colon[0,\, +\infty)\to [0,\, +\infty)$ belongs to class $\mathcal{K}$ if it is strictly increasing and $\alpha(0)=0$;  It further belongs to class $\mathcal{K}_\infty$ if it belongs to class $\mathcal{K}$ and  $\lim_{s\to \infty}\alpha(s)=\infty$.

\section{Preliminaries and Problem formulation}\label{sec:formulation}
In this section, we present preliminaries of partial stability and graph theory, and then the formulation of our problem.

\subsection{Partial stability}

To achieve optimal output consensus, we need to ensure the convergence of particular partial state 
of the closed-loop system rather than the full state. Such an issue is often referred to as partial stability (stabilization) \cite{vorotnikov1998partial}. Since the closed-loop system may have a continuum of equilibria, we introduce a modified version of partial stability as follows.  

Consider the nonlinear autonomous system
\begin{align}\label{sys:def-ps}
\dot{x}_1=f_1(x_1,\,x_2),\quad \dot{x}_2=f_2(x_1,\,x_2)
\end{align}  
where $x=\mbox{col}(x_1,\,x_2)$ with $x_1\in\R^{n_{x_1}}$, $x_2\in \R^{n_{x_2}}$ and the functions $f_1$, $f_2$ are sufficiently smooth.  Denote the equilibria set as $\mathcal{D}\triangleq \{x\mid f_1(x_1,\, x_2)={\bm 0},\,f_2(x_1,\, x_2)={\bm 0} \}$.

\begin{defn}\label{def:pss}
	System \eqref{sys:def-ps} is Lyapunov semistable with respect to $x_1$ (or briefly, $x_1$-semistable) at $x_1^{\star}$ if, for every $\e>0$, there exist $x_2^{\star}$  and $\delta>0$ such that $x^{\star}=\mbox{col}(x_1^{\star},\,{x_2^{\star}})\in \mathcal{D}$ and $||x(0)-x^{\star}||<\delta$ implies $||x_1(t)-x_1^{\star}||\leq \e$ for all $t\geq 0$. If for any $x(0)$, it further holds that $\lim_{t\to+\infty} ||x_1(t)-x_1^{\star}||=0$, this system is globally asymptotically $x_1$-semistable at $x_1^{\star}$.
\end{defn}

When $\mathcal{D}=\{{\bm 0}\}$, this definition is exactly the partial stability concept with respect to $x_1$ specified in [20, page~17].  The following lemma is slightly modified from Theorems 4.5 and 4.7 in \cite{haddad2008nonlinear} and its proof is omitted. 
\begin{lemj}\label{lem:ps+pss}
	Suppose that there exist a continuously differentiable function $V(x)$ and a constant vector $x_2^{\star}\in \R^{n_{x_2}}$ such that $x^{\star}=\mbox{col}(x_1^{\star},\,{x_2^{\star}})\in \mathcal{D}$, and along the trajectory of \eqref{sys:def-ps}, 
	\begin{align*}
	&\alpha(||x-x^{\star}||)\leq V(x)\leq \beta(||x-x^{\star}||) \\
	&\dot{V}(x)\leq -\gamma(||x_1-x^{\star}_1||)
	\end{align*} 
	for some functions $\alpha,\,\beta \in\mathcal{K}_\infty$ and $\gamma\in \mathcal{K}$. Then, system \eqref{sys:def-ps} admits well-defined bounded trajectories over $[0,\,+\infty)$ and is globally asymptotically $x_1$-semistable at $x_1^{\star}$. 
\end{lemj}

\subsection{Graph notion}
A weighted directed graph (digraph) is described by $\mathcal {G}=(\mathcal {N}, \mathcal {E}, \mathcal{A})$ with node set $\mathcal{N}=\{1,{\dots},N\}$ and edge set $\mathcal {E}$. $(i,\,j)\in \mathcal{E}$ denotes an edge from node $i$ to node $j$. The weighted adjacency matrix $\mathcal{A}=[a_{ij}]\in \mathbb{R}^{N\times N}$ is defined by $a_{ii}=0$ and $a_{ij}\geq 0$. Here $a_{ij}>0$ iff there is an edge $(j,\,i)$ in the digraph.  Node $i$'s neighbor set is defined as $\mathcal{N}_i=\{j\mid (j,\, i)\in \mathcal{E} \}$. We denote $\mathcal{N}_i^0=\mathcal{N}_i\cup \{i\}$.  A directed path is an alternating sequence $i_{1}e_{1}i_{2}e_{2}{\dots}e_{k-1}i_{k}$ of nodes $i_{l}$ and edges $e_{m}=(i_{m},i_{m+1}) \in\mathcal {E}$ for $l=1,2,{\dots},k$.  If there is a directed path between any two nodes, then the digraph is said to be strongly connected.  The in-degree and out-degree of node $i$ are defined by $d^{\mbox{in}}_i=\sum\nolimits_{j=1}^N a_{ij}$ and $d^{\mbox{out}}_i=\sum\nolimits_{j=1}^N a_{ji}$. A digraph is weight-balanced if $d^{\mbox{in}}_i=d^{\mbox{out}}_i$ for any $i\in \mathcal{N}$. The Laplacian of $\mathcal{G}$ is defined as $L\triangleq D^{\mbox{in}}-\mathcal{A}$ with $D^{\mbox{in}}=\mbox{diag}(d^{\mbox{in}}_1,\,\dots,\,d^{\mbox{in}}_N)$.  Note that $L{\bm 1}_N={\bm 0}_N$ for any digraph. If this digraph is weight-balanced, we have ${\bm 1}_N^\top L={\bm 0}_N^\top$ and  matrix  $\mbox{Sym}(L)\triangleq \frac{L+L^\top}{2}$ is positive semidefinite.  For a  weight-balanced and strongly connected digraph, we can order the eigenvalues of $\mbox{Sym}(L)$ as $0=\lambda_1<\lambda_2\leq \dots\leq \lambda_N$ and  have $ \lambda_2 I_{N-1}\leq M_2^\top \mbox{Sym}(L)M_2\leq \lambda_N I_{N-1}$. See \cite{godsil2001algebraic} for more details. 

\subsection{Problem formulation}

Consider a group of nonlinear systems modeled by  
\begin{align}\label{sys:agent}
\begin{split}
\dot{z}_i=\;&h_{i}(z_i,\,y_{i},\,w)\\
\dot{x}_{i}=\;&A_ix_i+B_i[g_{i}(z_i,\,x_{i},\,w)+ b_i(w) u_i]\\ 
y_i=\;&C_ix_{i},\quad i=1,\,\dotsm,\,N
\end{split}
\end{align}
where $\mbox{col}(z_i,\,x_i)$
is the state with $x_i=\mbox{col}(x_{i1},\,\dots,\,x_{i n_i}) \in \R^{n_i}$ and $z_i\in \R^{m_i}$, $u_i\in \R$ is the input, $y_i\in \R$ is the output, and $w\in \W\subset \R^{n_w}$ with $\W$ being compact and containing the origin. The triplet $(C_i,\,A_i,\,B_i)$ represents a chain of $n_{i}$ integrators in canonical form, that is,
\begin{align*}
A_i=\left[\begin{array}{c|c}
{\bm 0}_{n_i-1} &I_{n_{i}-1}\\\hline
0 & {\bm 0}_{n_i-1}^\top
\end{array}\right],\, B_i=\begin{bmatrix}
{\bm 0}_{n_i-1}\\
1
\end{bmatrix},\, 
C_i=\begin{bmatrix}
1\\
{\bm 0}_{n_i-1}
\end{bmatrix}^\top  %
\end{align*} 
Here $w$ and $z_i$ represent static and dynamic uncertainties of agent $i$, respectively. Different from \cite{tang2018acc}, the compact set $\W$  containing the static uncertainties is not necessarily known here.  It is assumed that all functions are sufficiently smooth and satisfy $h_{i}({\bm 0},\,0,\,w)=0$, $g_i({\bm 0},\, {\bm 0},\,w)=0$, $b_i(w)\geq b_0>0$ for all $w\in \W$ with some constant $b_0$.

We endow each agent output with a local cost function $f_i\colon\R\to\R$, and define the global cost function as the sum of all local costs, i.e., $f(y)=\sum_{i=1}^{N} f_i(y)$. For multi-agent system \eqref{sys:agent}, we aim to develop an algorithm such that all agent outputs achieve a consensus on the minimizer to this global cost function in a distributed fashion.  For this purpose, a digraph $\mathcal{G}=(\mathcal{N},\, \mathcal{E}, \,\mathcal{A})$  is used to describe the information communication relationships among agents with node set $\mathcal{N}=\{1,\,\dots,\, N\}$, edge set $\mathcal{E}\subset \mathcal{N}\times \mathcal{N}$, and weighted matrix $\mathcal{A}\in \R^{N\times N}$. An edge $(j,\,i )\in \mathcal{E}$ with weight $a_{ij}>0$ means that agent $i$ can get the information of agent $j$.
 
The considered distributed controller is described by 
\begin{align}\label{ctr:nominal}
\begin{split}
u_i=\;&\Xi_{i1}(\nabla f_i,\,x_j,\, \chi_j,\,j\in \mathcal{N}_i^0)\\
\dot{\chi}_i=\;&\Xi_{i2}(\nabla f_i,\,x_j,\,\chi_j,\,j\in \mathcal{N}_i^0)
\end{split}  
\end{align}
where $\chi_i\in \R^{q_i}$ is the compensator state and $\Xi_{i1}$, $\Xi_{i2}$ are smooth functions to be specified later. With these preparations, we formulate our problem explicitly as follows.
\begin{prop}\label{prop}
	For multi-agent system \eqref{sys:agent}, function $f_i$, digraph $\mathcal{G}$, and compact set $\W$,  find a controller of the form \eqref{ctr:nominal} such that, for each $w\in \W$ and each initial condition $\mbox{col}(z_i(0),\,x_i(0),\,{\chi}_i(0))\in \R^{m_i+n_i+q_i}$,
	\begin{itemize}
		\item[a)] the trajectory of the closed-loop system composed of \eqref{sys:agent} and \eqref{ctr:nominal} exists and is bounded over $[0,\,+\infty)$; 
		\item[b)] the outputs of agents satisfy $\lim_{t\to +\infty}|y_i(t)-y^{\star}|=0$ with $y^{\star}$ being optimal solution of 
		\begin{align}\label{opt:main}
		\min_{ y \in \R} \; f(y)=\sum\nolimits_{i=1}^{N} f_i(y)
		\end{align}
	\end{itemize}	
\end{prop}

\begin{remj}
	Compared to existing output consensus results \cite{kim2011output,liu2013distributed,rezaee2015average}, this problem further requires the outputs of agents to reach an agreement on the optimal point $y^{\star}$ specified by minimizing a cost function. In this sense, we say these agents achieving an optimal output consensus as in \cite{shi2013reaching,xie2019global,qiu2019distributed}.
 \end{remj}

This problem for single integrators has been coined as distributed optimization and investigated for many years. For high-order nonlinear agents, it is certainly  more challenging to achieve such an optimal output consensus, while the static and dynamic uncertainties bring extra technical difficulties in resolving this problem.

\section{Two-step Design Scheme}\label{sec:two}
In this section, we convert the optimal output consensus problem into a distributed partial stabilization problem by constructing optimal signal generators, giving rise to a two-step design scheme for solving Problem \ref{prop}.  

To begin with, several standing assumptions are listed.

\begin{assj}\label{ass:graph}
The digraph $\mathcal{G}$ is weight-balanced and strongly connected. 
\end{assj}

\begin{assj}\label{ass:convexity-strong}
	For each $i \in \mathcal{N}$, the function $f_i$ is twice continuously differentiable and satisfies that $\underline{l}_i  \leq \nabla^2 f_i(s) \leq \bar l_i$ with constants $0 < \underline{l}_i  \leq \bar l_i <+\infty$ for all $s\in \R$.
\end{assj}
\begin{assj}\label{ass:re-solution}
	For each $i\in  \mathcal{N}$, there exists a smooth function ${z}^{\star}_i(s,\,w)$ satisfying ${z}^{\star}_i(0,\,w)=0$ and $h_i({z}^{\star}_i(s,\,w),\,s,\, w)= 0$ for all $s \in\R$ and $w\in \R^{n_w}$.
\end{assj}

Assumption \ref{ass:graph} guarantees that each agent's information can be reached by any other agent. Assumption \ref{ass:convexity-strong} implies the existence and uniqueness of optimal solution to problem \eqref{opt:main} \cite{ruszczynski2011nonlinear}. Assumption \ref{ass:re-solution} can be interpreted as the solvability of regulator equations in the context of output regulation \cite{huang2004nonlinear}. These assumptions have been widely used in (distributed) coordination for multi-agent systems \cite{jakovetic2015linear, deng2018ijrnc, wang2016distributed,zhang2017distributed, su2015cooperative}.

Consider an optimal consensus problem for a group of single integrators with the same optimal requirement \eqref{opt:main}
\begin{align}
\dot{r}_i=\mu_i
\end{align}
If this auxiliary problem is solved by some chosen $\mu_i$, we only need to drive agent $i$ to track the generated signal $r_i(t)$ to achieve the optimal output consensus for agent \eqref{sys:agent}. 

Since the Laplacian $L$ of digraph $\mathcal{G}$ is asymmetric, the generator in \cite{tang2017cyb} fails to reproduce $y^{\star}$ without the information of $L^\top$. Motivated by \cite{gharesifard2014distributed}, we present a candidate of optimal signal generator for problem \eqref{opt:main} as follows
	\begin{align}\label{sys:generator}
	\mu_i=\;&-\alpha \nabla f_i(r_i)-\beta \sum\nolimits_{j=1}^{N}a_{ij}(r_i-r_j)-\sum\nolimits_{j=1}^{N}a_{ij}(v_i-v_j)\nonumber\\
	\dot{v}_i=\;&\alpha \beta  \sum\nolimits_{j=1}^{N}a_{ij}(r_i-r_j)
	\end{align}
where $\alpha,\,\beta$ are constants to be specified later.  Putting it into a compact form gives
\begin{align}\label{sys:composite-osg}
\begin{split}
\dot{r}=\;&-\alpha \nabla \tilde f(r)- \beta Lr-Lv,\quad \dot{v}=\alpha \beta Lr
\end{split}
\end{align}
where $r=\mbox{col}(r_1,\,\dots,\,r_N)$, $v=\mbox{col}(v_1,\,\dots,\,v_N)$, and function $\tilde f(r)\triangleq \sum\nolimits_{i=1}^Nf_i(r_i)$ is $\underline{l}$-strongly convex while its gradient $ \nabla\tilde f(r)$ is $\bar l$-Lipschitz with $\bar l=\max_i\{\bar l_i\}$ and $\underline{l}=\min_i\{\underline{l}_i\}$. 

Let $\mbox{col}(r^{\star},\,v^{\star})$ be the equilibrium point of system \eqref{sys:composite-osg}. It is verified that $r^{\star}={\bm 1}_N y^{\star}$ under Assumptions \ref{ass:graph} and \ref{ass:convexity-strong} by Theorem 3.27 in \cite{ruszczynski2011nonlinear}. For \eqref{sys:composite-osg} , we have the following interesting result. 

\begin{lemj}\label{lem:generator}
	Suppose Assumptions \ref{ass:graph}--\ref{ass:convexity-strong} hold and let $\alpha\geq \max\{1,\,\frac{1}{\underline{l}},\,\frac{2\bar l^2}{\underline{l}\lm_2}\}$, $\beta\geq \max\{1,\, \frac{1}{\lambda_2},\,\frac{6\alpha^2\lambda_N^2}{\lambda_2^2} \}$. 	Then, system \eqref{sys:composite-osg} admits well-defined bounded trajectories over $[0,\,+\infty)$ and is globally asymptotically $r$-semistable at ${\bm 1}_N y^{\star}$. Moreover, $r_i(t)$ approaches  $y^{\star}$ exponentially as $t\to +\infty$ for $i\in\mathcal{N}$.
\end{lemj}
\pb Briefly, we utilize Lemma \ref{lem:ps+pss} to complete the proof. Let $M_L= M_2^\top LM_2$ and $v^{\star}= -\alpha M_2 M_L^{-1}M_2^\top \nabla \tilde f(r^*)$. It can be verified that $\mbox{col}(r^{\star},\,v^{\star})$ is an equilibrium of system \eqref{sys:composite-osg}. 
 
Perform the coordinate transformation: $\bar r_1=M_1^\top (r-r^\star)$, $\bar r_2=M_2^\top (r-r^\star)$, $\bar v_1=M_1^\top (v-v^\star)$, and $\bar v_2=M_2^\top[( v+\alpha r)-( v^\star+\alpha r^\star)]$. It follows that $\dot{\bar v}_1= 0$ and  
\begin{align}\label{sys:composite-osg-reduced}
\begin{split}
\dot{\bar r}_1=\;&-\alpha M_1^\top {\bm \Pi}\\
\dot{\bar r}_2=\;&-\alpha M_2^\top {\bm \Pi}-\beta M_L \bar r_2 + \alpha M_L \bar r_2-M_L\bar v_2\\
\dot{\bar v}_2=\;&-\alpha M_L {\bar v}_2+\alpha^2 M_L \bar r_2-\alpha^2 M_2^\top {\bm \Pi} 
\end{split}
\end{align}
where  ${\bm \Pi}\triangleq \nabla \tilde f(r)-\nabla \tilde f(r^\star)$.  Let $\bar r=\mbox{col}(\bar r_1,\,\bar r_2)$, and $V_{\rm o}(r,\, v)=\bar r^\top \bar r + \frac{1}{\alpha^3} \bar v_1^\top \bar v_1+\frac{1}{\alpha^3}\bar v_2^\top \bar v_2$ in this new coordinate with $\alpha>0$ to be specified later. The first inequality in Lemma \ref{lem:ps+pss} apparently hold. On the other hand, by Young's inequality, the time derivative of $V_{\rm o}$ along the trajectory of \eqref{sys:composite-osg} satisfies 
\begin{align*}
\dot{V}_{\rm o}
=\;&-2\alpha (r-r^\star)^\top {\bm \Pi}+ 2\bar r_2^\top [-\beta M_L \bar r_2 + \alpha M_L \bar r_2-M_L\bar v_2 ]\\
&+\frac{2}{\alpha^3}\bar v_2^\top[-\alpha M_L {\bar v}_2+\alpha^2 M_L \bar r_2-\alpha^2 M_2^\top {\bm \Pi} ]\\
\leq \;& -2\alpha\underline{l}||\bar r||^2-2\beta\lambda_2||\bar r_2||^2+2\alpha  \lambda_N  ||\bar r_2||^2+2 \lambda_N  ||\bar r_2||||\bar v_2||\\
&-\frac{2\lambda_2}{\alpha^2}||\bar v_2||^2+ \frac{2}{\alpha} \lambda_N ||\bar r_2|| ||\bar v_2||+\frac{2\bar l}{\alpha}||\bar v_2||||\bar r||\\
\leq\;& -(2\alpha\underline{l} -\frac{3\bar l^2}{\lambda_2} )||\bar r||^2-\frac{\lambda_2}{\alpha^2}||\bar v_2||^2\\
& - (2\beta\lambda_2-2\alpha\lambda_N-\frac{3\alpha^2\lambda_N^2}{\lambda_2}-\frac{3\lambda_N^2}{\lambda_2})||\bar r_2||^2\\
&\leq - \frac{1}{2}||\bar r||^2-\frac{1}{2\alpha^3}||\bar v_2||^2\triangleq W_{\rm o}(\bar r,\, \bar v_2)
\end{align*}
According to Lemma \ref{lem:ps+pss}, we conclude the boundedness of all trajectories over $[0,\,+\infty)$ and its $r$-semistability of system \eqref{sys:composite-osg} at ${\bm 1}_N y^{\star}$. By further considering the reduced-order system \eqref{sys:composite-osg-reduced} with a Lyapunov function $W_{\rm o}(\bar r,\,\bar v_2)$,  one can obtain that $\dot{W}_{\rm o}\leq -\frac{1}{2}W_{\rm o}$ along the trajectories of \eqref{sys:composite-osg-reduced}. Recalling Theorem 4.10 in \cite{khalil2002nonlinear},  $W_{\rm o}(\bar r(t),\,\bar v_2(t))$ and $\bar r(t)$ must exponentially converge to $0$ as $t$ goes to infinity. The proof is complete.
\pe

\begin{remj}
	The optimal signal generator \eqref{sys:composite-osg} is a modified version of the augmented Lagrangian method solving problem \eqref{opt:main} in \cite{gharesifard2014distributed}. Here we add an extra parameter $\alpha$ to simplify both the synthesis and its analysis. Compared with the results for digraphs in \cite{kia2015distributed,xie2019global,qiu2019distributed}, our algorithm is initialization-free to generate the optimal point $y^*$. This makes it possible to work in a scalable manner, which might be favorable for dynamic networks with leaving-off and plugging-in of agents. 
\end{remj}

\begin{remj}
	In our design, we use the knowledge of $\lambda_2$ and $\lambda_N$ as that in \cite{zhang2017distributed,deng2018ijrnc} to compensate the asymmetry of directed information flows. 
    It should be mentioned that these values can be computed by existing algorithms beforehand, e.g.,\cite{charalambous2015distributed}. 
\end{remj}

Under Assumption \ref{ass:re-solution}, we denote ${x}^{\star}_i(r_i)=\mbox{col}(r_i,\,{\bm 0}_{n_i-1})$, ${u}^{\star}_i(r_i,\,w)=-\frac{g_i({z}^{\star}_i(r_i,\,w),\,{x}^{\star}_i(r_i),\,w)}{b_i(w)}$ and perform the coordinate transformation: $\bar z_i=z_i-{z}^{\star}_i(r_i,\,w)$,\, $\bar x_{i}=x_{i}-{x}^{\star}_i(r_i)$. This leads to an interconnected error system as follows
\begin{align}\label{sys:agent-error}
\dot{\bar z}_i=\;&\bar h_{i}(\bar z_i,e_{i},r_i,w)-\frac{\partial {z}^{\star}_i}{\partial r_i} \mu_i  \nonumber\\
\dot{\bar x}_{i}=\;&A_i\bar x_i+B_i[\bar g_{i}(\bar z_i,\,\bar x_{i},\,r_i,\,w)  \nonumber \\
&+b_i(w)(u_i- {u}^{\star}_i(r_i,\,w))]- E_i \mu_i  \\ 
e_i=\;&C_i\bar x_{i},\quad i\in\mathcal{N}\nonumber
\end{align}
where $E_i=\mbox{col}(1,\,{\bm 0}_{n_i-1})$ and 
\begin{align*}
\bar h_{i}(\bar z_i,\, e_{i},\, r_i,\,w)=\;&h_{i}(z_i,\, y_i,\, w)-h_i({z}^{\star}_i(r_i,\,w),\, r_i,\,w)\\
\bar g_i(\bar z_i,\, \bar x_i,\, r_i,\, w)=\;&g_{i}(z_i,\, x_i,\,w)-g_{i}({z}^{\star}_i(r_i,\,w),\,{x}^{\star}_i(r_i),\,w)
\end{align*}
It can be verified that $\bar h_{i}({\bm 0},\, 0,\, r_i,\,w)=0, \, \bar g_i({\bm 0},\, {\bm 0},\, r_i,\, w)=0$ for all $r_i\in \R$ and $w\in \R^{n_w}$. 

Attaching the optimal signal generator \eqref{sys:composite-osg} to error system \eqref{sys:agent-error} yields an augmented system associated with Problem \ref{prop}. A key lemma is obtained to assist us in solving the optimal output consensus problem.
\begin{lemj}\label{lem:conversion}	
Suppose Assumptions \ref{ass:graph}--\ref{ass:re-solution} hold and there exists a smooth controller of the form 
\begin{align}\label{ctr:nominal-partial}
\begin{split}
u_i=\;&\Xi_{i1}^{\rm o}(\bar x_j,\,r_j,\,\chi_j^{\rm o},\,j\in \mathcal{N}_i^0)\\
\dot{\chi}_i^{\rm o}=\;&\Xi_{i2}^{\rm o}(\bar x_j, \,r_j,\,\chi_j^{\rm o},\, j\in \mathcal{N}_i^0)
\end{split}
\end{align}
solving the distributed partial stabilization problem of the augmented system composed of \eqref{sys:composite-osg} and \eqref{sys:agent-error} in the sense that the closed-loop system composed of \eqref{sys:composite-osg}, \eqref{sys:agent-error}, and \eqref{ctr:nominal-partial} admits well-defined bounded trajectories over $[0,\,+\infty)$  and is globally
asymptotically $e_i$-semistable at $0$.  Then,  Problem \ref{prop} can be solved by a controller composed of \eqref{sys:generator} and \eqref{ctr:nominal-partial}.
\end{lemj}
\pb  
Under the lemma condition, we can confirm that trajectories of all agents are well-defined bounded over $[0,\, +\infty)$ and $\lim_{t\to +\infty} e_i(t)=0$ for any initial condition $\mbox{col}(\bar z_i(0),\, \bar x_i(0),\, {\chi}_i^{\rm o}(0),\,r(0),\,v(0))$.  Note that  $|y_i(t)-y^{\star}|\leq|e_i(t)|+|r_i(t)-y^{\star}|$ by the triangle inequality. This together  with Lemma 2 ensures that $\lim_{t\to +\infty}|y_i(t)-y^{\star}|=0$. 
\pe

\begin{remj}
	Based on Lemma \ref{lem:conversion}, our optimal output consensus problem for multi-agent system \eqref{sys:agent} is converted into a distributed partial stabilization problem of certain interconnected augmented systems. 	As the considered nonlinear multi-agent system \eqref{sys:agent} is further subject to static and dynamic uncertainties, the associated partial stabilization design is more challenging than relevant results obtained in \cite{qiu2019distributed, tang2018ijrnc, wang2016distributed, tang2017cyb}. On the other hand, existing designs presented in \cite{vorotnikov1998partial,haddad2008nonlinear} are not applicable for such complicated uncertainties and the partial stabilization problem itself is nontrivial even for a single nonlinear system. Thus, we have to seek a robust distributed partial stabilization design method for the augmented systems.
\end{remj}

\section{Main Result}\label{sec:main}

In this section, we focus on the subsequent partial stabilization problem of the  augmented system composed of  \eqref{sys:composite-osg} and \eqref{sys:agent-error} and eventually solve the optimal output consensus problem for multi-agent system \eqref{sys:agent}.

To this end, we make an extra assumption imposing a mild minimum-phase condition widely used in nonlinear stabilization problems \cite{xu2010robust,tang2015distributed,su2015cooperative}.
\begin{assj}\label{ass:zero-dynamics}
	For each $ i\in \mathcal{N}$, there exists a continuously differentiable function $W_{i\bar z}(\bar z_i)$ such that, for all $r_i\in \R$ and $w\in \W$, along the trajectory of system \eqref{sys:agent-error}, 
	\begin{align}\label{eq:zero-dynamics}
	\begin{split}
	&\underline{\alpha}_{i}(||{\bar z}_i||)\leq W_{i\bar z}({\bar z}_i)\leq {\bar \alpha}_i(||{\bar z}_i||) \\ 
	&\dot{W}_{i\bar z}\leq -\alpha_i(||\bar z_i||)+ \sigma_{ie} \gamma_{ie}(e_i) e_i^2+\sigma_{i\mu} \gamma_{i\mu}(r_i)\mu_i^2
	\end{split}
	\end{align}
	for some known smooth functions $\underline{\alpha}_i$, ${\bar \alpha}_i$, $\alpha_i \in \mathcal{K}_\infty$, $\gamma_{ie}$, $\gamma_{ir}>1$, and  unknown constants $\sigma_{ie}$, $\sigma_{i\mu}>1$ with $\alpha_i$ satisfying $\limsup_{s \to 0+}\frac{\alpha_i^{-1}(s^2)}{s}<+\infty$.
\end{assj}

Due to the presence of uncertain parameter $w$, the feedforward term ${u}^{\star}_i(r_i,\,w)$ is unavailable for feedback.  To tackle this issue, we introduce a dynamic compensator as follows
\begin{align*}
\dot{\eta}_i=-\kappa_i(r_i)\eta_i+u_i
\end{align*}
where $\kappa_i(r_i)>0$ is a smooth function to be specified later. Here, $\kappa_i(r_i)$ is a scaling factor to handle nonlinear functions of $r_i$. This compensator reduces to an internal model when $\kappa_i(r_i)$ is constant \cite{huang2004nonlinear}. 

Consider the error system \eqref{sys:agent-error}. For $n_i\geq 2$, choose constants $k_{ij}$ such that the polynomial $p_i(\lambda)=\sum_{j=1}^{n_i-1}k_{ij} \lambda^{j-1}+\lambda^{n_i-1}$ is Hurwitz. Let  $\xi_i=\mbox{col}(\bar x_{i1},\,\dots,\,\bar x_{i n_i-1})$, $\zeta_i=\sum_{j=1}^{n_i-1}k_{ij} \bar x_{ij}+\bar x_{in_i}$, and $\beta_i(\eta_i,\,r_i) \triangleq \kappa_i(r_i)\eta_i$. Performing coordinate and input transformations: $\bar \eta_i=\eta_i-\frac{{u}^{\star}_i(r_i,\,w)}{\kappa_i(r_i)}-b_i^{-1}(w)\zeta_i$ and $\bar u_i=u_i-\beta_i(\eta_i,\,r_i)$ gives a composite system in the following form
\begin{align}\label{sys:agent-error:translated-im}
\begin{split}
\dot{\bar z}_i=\;& \bar h_{i}(\bar z_i,\,e_{i},\,r_i,\,w)-\frac{\partial {z}^{\star}_i}{\partial r_i} \mu_i\\
\dot{\xi}_i=\;&  A_i^{\rm o} {\xi}_i+ B_i^{\rm o}\zeta_i- E_i^{\rm o} \mu_i\\
\dot{\bar \eta}_i=\;& -\kappa_i(r_i)\bar \eta_i+\tilde g_i(\bar z_i,\,{\xi}_i,\, \zeta_i,\, r_i,\, w)+\psi_i(r_i,\,w)\mu_i\\
\dot{\zeta}_i=\;&\check{g}_i(\bar z_i,\,\xi_i,\,\bar \eta_i,\,\zeta_i,\, r_i,\,w)+b_i(w)\bar u_i-k_{i1}\mu_i
\end{split}
\end{align}
where 
\begin{align*}
A_i^{\rm o}=\;&\left[\begin{array}{c|c}
{\bm 0}_{n_i-2}&I_{n_i-2}\\\hline
-k_{i1}&-k_{i2},\dots,-k_{in_i-1}
\end{array}\right]\\
B_i^{\rm o}=\;&\begin{bmatrix}
{\bm 0}_{n_i-2}\\
1
\end{bmatrix},~~E_i^{\rm o}=\begin{bmatrix}
1\\
{\bm 0}_{n_i-2}
\end{bmatrix}\\
\tilde g_i=\;&-\frac{1}{b_i(w)}[\hat g_i(\bar z_i,{\xi}_i,\zeta_i,r_i,w)+\kappa_i(r_i)\zeta_i]\\
\psi_i=\;&\frac{{u}^{\star}_i(r_i,\, w)}{\kappa^2_i(r_i)}\frac{\partial {\kappa}_i(r_i)}{\partial r_i}-\frac{1}{\kappa_i(r_i)}\frac{\partial {u}^{\star}_i(r_i,\, w)}{\partial r_i}+\frac{k_{i1}}{b_i(w)}\\
\check{g}_i=\;&\kappa_i(r_i)\zeta_i+b_i(w)\kappa_i(r_i)\bar \eta_i+\hat g_i(\bar z_i,\xi_i,\zeta_i,r_i,w)\\
\hat g_i=\;&-k_{in_i-1}k_{i1}\bar x_{i1}+\sum\nolimits_{j=2}^{n_i-1} (k_{ij-1}-k_{in_i-1}k_{ij})\bar x_{ij}\\
&+k_{in_i-1}\zeta_i+\bar g_{i}(\bar z_i,\bar x_{i},r_i,w)
\end{align*}
It can be verified that $\hat g_i({\bm 0},\, {\bm 0},\, 0,\, r_i,\, w)=0$, $\tilde g_i({\bm 0},\,{\bm 0},\,0,\, r_i,\, w)=0$, and $\check{g}_i({\bm 0},\,{\bm 0},\,0,\,0,\, r_i,\,w)=0$ for all $r_i\in \R$ and $w\in \R^{n_w}$. Denote $\tilde z_i=\mbox{col}(\bar z_i,\,\xi_i)$ and $\hat z_i=\mbox{col}(\tilde z_i,\,\bar \eta_i)$. For $n_i=1$, the ${\xi}_{i}$-subsystem vanishes and we let $\tilde z_i=\bar z_i$, $\zeta_i=\bar x_{i1}$ for consistency. 

According to Lemma 11.1(iv) in \cite{chen2015stabilization} and by completing the square, there exist some known  smooth functions $\hat\phi^0_{i1},\,\hat \phi_{i2},\,\hat \phi_{i3}>1$ such that, for all $r_i\in \R$ and $w\in \W$, 
\begin{align}\label{eq:growth11}
||\hat g_i(\tilde z_i,\, \zeta_i,\, r_i,\, w)||^2\leq \hat \phi_{i1}^0(r_i,\,w)[\hat \phi_{i2}(\tilde  z_i)||\tilde  z_i||^2+\hat \phi_{i3}(\zeta_i)\zeta_i^2]
\end{align}
By Lemma 11.1(i) in \cite{chen2015stabilization}, there exist some known smooth functions $\hat \phi_{i1},\,\hat \phi_{i4}>1$ and unknown constants $\hat c_{ig}$,\, $\hat\ell_{i\psi}>1$ satisfying 
\begin{align}\label{eq:growth12}
\hat \phi^{0}_{i1}(r_i,\,w)\leq \hat c_{ig} \hat \phi_{i1}(r_i),\quad \psi_i^2(r_i,\,w)\leq \hat\ell_{i\psi}\hat \phi_{i4}(r_i)
\end{align}
It follows that, for all $r_i\in \R$ and $w\in \W$, 
\begin{align}\label{eq:growth1}
||\hat g_i(\tilde z_i,\,\zeta_i,\,r_i,\,w)||^2\leq \hat c_{ig} \hat \phi_{i1}(r_i)[\hat \phi_{i2}(\tilde z_i)||\tilde z_i||^2+\hat \phi_{i3}(\zeta_i)\zeta_i^2] 
\end{align}  
Similarly, one can determine some known smooth functions $\check \phi_{i1},\,\check \phi_{i2},\,\check \phi_{i3}>1$ and unknown constant $\check c_{ig}>1$ such that, for all $r_i\in \R$ and $w\in \W$, 
\begin{align}\label{eq:growth2}
||\check{g}_i(\hat z_i,\,\zeta_i,\,r_i,\, w)||^2 \leq \check c_{ig} \check \phi_{i1}(r_i)[\check \phi_{i2}(\hat  z_i)||\hat z_i||^2+\check \phi_{i3}(\zeta_i)\zeta_i^2] 
\end{align}
We claim the $\hat z_i$-subsystem admits the following property. 
\begin{lemj}\label{lem:zero-dynamics}
	For each $i\in \mathcal{N}$, let  $\kappa_i(r_i)\geq \hat \phi_{i1}(r_i)+1$. Then, there exists a continuously differentiable function $W_i(\hat z_i)$ such that, for all $r_i\in \R$ and $w\in \W$, along the trajectory of \eqref{sys:agent-error:translated-im},
	\begin{align*}
	&\hat {\underline \alpha}_{i}(||\hat z_i||) \leq W_i(\hat z_i)\leq \hat {\bar \alpha}_i(||\hat z_i||) \\ 
	&\dot{W}_i(\hat z_i)\leq -||\hat z_i||^2 +\hat \sigma_{i\zeta}\hat \gamma_{i\zeta}(\zeta_i,\,r_i)\zeta_i^2+\hat \sigma_{i\mu}\hat \gamma_{i\mu}(\mu_i,\,r_i)\mu_i^2\nonumber
	\end{align*}
	for some known smooth functions $\hat {\underline{\alpha}}_i,\,\hat {\bar \alpha}_i \in \mathcal{K}_\infty$, $\hat \gamma_{i\zeta},\,\hat \gamma_{i\mu}>1$, and unknown constants $\hat \sigma_{i\zeta}$, $\hat \sigma_{i\mu}>1$.
\end{lemj}

The proof of Lemma \ref{lem:zero-dynamics} is put in {\bf Appendix}.

Motivated by \cite{xu2010robust,tang2015distributed}, we let $\bar u_i=-\theta_i\rho_i(\zeta_i,\,r_i)\zeta_i$ with $\dot{\theta}_{i}=\tau_i(\zeta_i,\,r_i)$. Here, $\rho_i$ and $\tau_i$ are positive smooth functions to be specified later and $\theta_i$ is a dynamic gain to handle the unknown boundaries of  static uncertainties. For simplicity, we set $\theta_i(0)=0$.  The developed partial stabilizer for the augmented system  \eqref{sys:composite-osg}--\eqref{sys:agent-error} is consequently
\begin{align}\label{ctr:adaptive-partial}
\begin{split}
u_i=\;&-\theta_i\rho_i(\zeta_i,\,r_i)\zeta_i+\kappa_i(r_i)\eta_i\\
\dot{\eta}_i=\;&-\kappa_i(r_i)\eta_i+u_i\\
\dot{\theta}_{i}=\;&\tau_i(\zeta_i,\,r_i)
\end{split}
\end{align}
It is of the form \eqref{ctr:nominal-partial} and distributed in the sense of using each agent's own and neighboring information. 

We are ready to present our main theorem. 
\begin{thmj}\label{thm:main}
	Under Assumptions \ref{ass:graph}--\ref{ass:zero-dynamics}, there exist positive constants $\alpha,\,\beta$ and smooth functions $\kappa_{i}(r_i)$, $\rho_i(\zeta_i,\,r_i)$, $\tau_i(\zeta_i,\,r_i)$ such that Problem \ref{prop} for multi-agent system \eqref{sys:agent} is solved by a distributed controller of the following form
	\begin{align}\label{ctr:adaptive}
	u_i=\;&-\theta_i\rho_i(\zeta_i,\,r_i)\zeta_i+\kappa_i(r_i)\eta_i \nonumber\\
	\dot{\eta}_i=\;&-\kappa_i(r_i)\eta_i+u_i \nonumber\\
	\dot{\theta}_{i}=\;&\tau_i(\zeta_i,\,r_i)\\
	\dot{r}_i=\;&-\alpha \nabla f_i(r_i)-\beta \sum\nolimits_{j=1}^{N}a_{ij}(r_i-r_j)-\sum\nolimits_{j=1}^{N}a_{ij}(v_i-v_j)\nonumber\\
	\dot{v}_i=\;&\alpha \beta  \sum\nolimits_{j=1}^{N}a_{ij}(r_i-r_j) \nonumber
	\end{align}
\end{thmj}
\pb Set $\alpha,\,\beta$ and $\kappa_{i}(r_i)$ as in Lemmas \ref{lem:generator} and \ref{lem:zero-dynamics}. By Lemma \ref{lem:conversion}, we are left to show the following closed-loop system admits well-defined bounded trajectories for $t\geq 0$  and is globally asymptotically $e_i$-semistable at $0$.
\begin{align}\label{sys:agent-error-proof}
\dot{\hat z}_i=\;& \hat h_{i}(\hat z_i,\, \zeta_{i},\, r_i, \, w,\, \mu_i)\nonumber\\
\dot{\zeta}_i=\;& \check{g}_i(\hat z_i,\,\zeta_i,\, r_i,\,w)-\theta_ib_i(w)\rho_i(\zeta_i,\,r_i)\zeta_i-k_{i1}\mu_i\nonumber\\
\dot{\theta}_{i}=\;&\tau_i(\zeta_i,\,r_i) \\
\dot{r}_i=\;&-\alpha \nabla f_i(r_i)-\beta \sum\nolimits_{j=1}^{N}a_{ij}(r_i-r_j)-\sum\nolimits_{j=1}^{N}a_{ij}(v_i-v_j)\nonumber\\
\dot{v}_i=\;&\alpha \beta  \sum\nolimits_{j=1}^{N}a_{ij}(r_i-r_j) \nonumber
\end{align}
where function $\hat h_i$ is determined by \eqref{sys:agent-error:translated-im} and we simply denote $\check{g}_i(\hat z_i,\,\zeta_i,\, r_i,\,w)\triangleq \check{g}_i(\bar z_i,\,\xi_i,\,\bar \eta_i,\,\zeta_i,\, r_i,\,w)$ to save notations. 

The proof is divided into two steps. 

{\it Step 1}: we consider the first three subsystems and seek certain disturbance attenuation performance with  $\mu_i$ as its disturbance by choosing $\rho_i$ and $\tau_i$. 

First, by Lemma \ref{lem:zero-dynamics}, we apply the changing supply functions technique \cite{sontag1995changing} and conclude that, for any given smooth function $\hat \Delta_i(\hat z_i)>0$, there exists a continuously differentiable function $W^1_{i}(\hat z_i)$ such that, along the trajectory of \eqref{sys:agent-error-proof},
\begin{align*}
&\hat {\underline \alpha}^1_{i}(||\hat   z_i||) \leq W_i^1(\hat z_i)\leq \hat {\bar \alpha}^1_i(||\hat  z_i||)\\
&\dot{W}^1_i\leq -\hat \Delta_i(\hat z_i)||\hat  z_i||^2+ \hat \sigma^1_{i\zeta} \hat \gamma^1_{i\zeta}(\zeta_i,r_i)\zeta_i^2+ \hat \sigma^1_{i\mu} \hat \gamma^1_{i\mu}(\mu_i,r_i)\mu_i^2
\end{align*}
for some known smooth functions $\hat {\underline\alpha}^1_i,\,\hat {\bar \alpha}^1_i \in \mathcal{K}_\infty$, $\hat \gamma^1_{i\zeta},\,\hat \gamma^1_{i\mu}>1$, and unknown constants $\hat \sigma^1_{i\zeta},\,\hat \sigma^1_{i\mu}>1$.

Second, let $V_i(\hat z_i,\,\zeta_i,\bar \theta_i)=\hat \ell_i  W_i^1(\hat z_i)+\zeta_i^2+\bar \theta_i^2$, where $\bar \theta_i=\theta_i- \Theta_i$ with $\Theta_i,\,\hat \ell_i>0$ to be specified later. It is positive definite and radially unbounded, and moreover satisfies
\begin{align*}
\dot{V}_i\leq\;&- \hat \ell_i [\hat \Delta_i(\hat z_i)||\hat  z_i||^2-\hat \sigma^1_{i\zeta} \hat \gamma^1_{i\zeta}(\zeta_i,r_i)\zeta_i^2-\hat \sigma^1_{i\mu} \hat \gamma^1_{i\mu}(\mu_i,r_i)\mu_i^2]\\
&+2\zeta_i[\check g_i(\hat z_i,\,\zeta_i,\, r_i,\, w)-\theta_ib_i(w)\rho_i(\zeta_i,\,r_i)\zeta_i-k_{i1}\mu_i ]\\
&+2(\theta_i-\Theta_i)\tau_i(\zeta_i,\, r_i)
\end{align*}
Recalling inequality \eqref{eq:growth2}, we complete the square and have  
\begin{align*}
\dot{V}_i\leq\;& -[ \hat \ell_i \hat \Delta_i(\hat z_i)-\check c_{ig} \check \phi_{i2}(\hat  z_i)]||\hat  z_i||^2-[2\theta_i b_i(w)\rho_i(\zeta_i,\,r_i)\\
&-\check \phi_{i1}(r_i)-\hat \ell_i \hat \sigma^1_{i\zeta}  \hat \gamma^1_{i\zeta}(\zeta_i,\,r_i)- \check c_{ig}\check \phi_{i3}(\zeta_i)-1]\zeta_i^2\\
&+ [\hat \ell_i \hat \sigma^1_{i\mu}\hat \gamma^1_{i\mu}(\mu_i,\,r_i)+k^2_{i1}]\mu_i^2+2(\theta_i-\Theta_i)\tau_i(\zeta_i,\, r_i)
\end{align*}
Choosing 
\begin{align}\label{eq:control-design-thm1}
\begin{split}
&\hat \ell_i \geq \check c_{ig}, \quad \hat \Delta_i(\hat z_i)\geq \check \phi_{i2}(\hat  z_i)+1\\ 
&\rho_i(\zeta_i,\,r_i)  \geq\hat \gamma^1_{i\zeta}(\zeta_i,\,r_i)+\check \phi_{i1}(r_i)+\check \phi_{i3}(\zeta_i)+2\\
&\tau_i(\zeta_i,\,r_i)=\rho_i(\zeta_i,\,r_i)\zeta_i^2, \quad \Theta_i  \geq \frac{1}{2b_0} \max\{\hat \ell_i \hat \sigma^1_{i\zeta},\, \check c_{ig} \}
\end{split}
\end{align} 
gives  $\dot{V}_i\leq -||\hat z_i||^2-\zeta_i^2+ [\hat \ell_i \hat \sigma^1_{i\mu}\hat \gamma^1_{i\mu}(\mu_i,\,r_i)+k^2_{i1}]\mu_i^2$. By Lemma \ref{lem:generator} and the smoothness of $\hat \gamma^1_{i\mu}$, there exists a constant  $c_{i\mu}>0$ satisfying $\ell_i^1 \hat \sigma^1_{i\mu}\hat \gamma^1_{i\mu}(\mu_i,\,r_i)+k^2_{i1}\leq c_{i\mu}$,
which further implies  
\begin{align*}
\dot{V}_i\leq -||\hat z_i||^2-\zeta_i^2+c_{i\mu} \mu_i^2
\end{align*}

{\it Step 2}: we show that the closed-loop system \eqref{sys:agent-error-proof} admits well-defined bounded trajectories for $t>0$ and is globally asymptotically $e_i$-semistable at $0$. 
 
Note that the equilibria set of \eqref{sys:agent-error-proof} is specified by $\mathcal{D}=\{\mbox{col}(\hat z,\, \zeta,\, \theta,\, r,\, v)\mid\hat z={\bm 0},\,\zeta={\bm 0},\,r={\bm 1}_N y^{\star},\, v=v^{\star}+ l_v {\bm 1}_N\}$ with an arbitrary constant $l_v$. For $e_i=0$, we set $\Theta^{\star}=\mbox{col}(\Theta_1,\,\dots,\,\Theta_N)$ and verify that $\mbox{col}({\bm 0},\,{\bm 0},\,\Theta^{\star},\,{\bm 1}_N y^{\star},\,v^{\star})$ is an equilibrium of system \eqref{sys:agent-error-proof}. 

From the proof of Lemma \ref{lem:generator}, we know that the function $W_{\rm o}(\bar r,\,\bar v_2)$ defined thereof satisfies $\bar \ell_{1}||\mbox{col}(\bar r,\,\bar v_2)||^2\leq W_{\rm o}(\bar r,\,\bar v_2)\leq \bar \ell_{2} ||\mbox{col}(\bar r,\,\bar v_2)||^2$ and $\dot{W}_{\rm o} \leq - \bar \ell_{3} ||\mbox{col}(\bar r,\,\bar v_2)||^2$ for  some constants $\bar \ell_{1},\,\bar \ell_{2},\, \bar \ell_{3}>0$. Due to the Lipschitzness of ${\bm \Pi}$ in $\bar r$, $\mu_i$ is also Lipschitz in $\mbox{col}(\bar r,\,\bar v_2)$. Thus, there exists a constant $\bar \ell_4>0$ such that $\sum_{i=1}^N c_{i\mu} \mu_i^2\leq \bar \ell_4 \bar \ell_3||\mbox{col}(\bar r,\,\bar v_2)||^2$.  

Let $V=\sum_{i=1}^N V_i+\bar \ell_4 V_{\rm o}$ with $V_{\rm o}$ defined in the proof of Lemma \ref{lem:generator}. The first condition in Lemma \ref{lem:ps+pss} is verified.  Taking the time derivative of $V$ along the trajectory of  \eqref{sys:agent-error-proof} gives
\begin{align*}
\dot{ V}\leq\;& -||\hat z||^2-\zeta^2+\sum_{i=1}^N c_{i\mu} \mu_i^2  - \bar \ell_{4} \bar \ell_{3} ||\mbox{col}(\bar r,\,\bar v_2)||^2\\
\leq\;&  -||\hat z||^2-\zeta^2
\end{align*}
This implies the second inequality in Lemma \ref{lem:ps+pss}.  Overall, the function $V$ indeed satisfies the conditions in Lemma \ref{lem:ps+pss}. This guarantees the trajectory's boundedness over $[0,\,+\infty)$ and the global asymptotic $e_i$-semistability of system \eqref{sys:agent-error-proof} at ${0}$.  By Lemma \ref{lem:conversion}, we complete the proof.  
\pe 

	\begin{remj}\label{rem:control}
		The developed optimal consensus control \eqref{ctr:adaptive} is of a high-gain type to handle the uncertainties. The parameters and functions can be sequentially constructed. Firs, we choose $\alpha$,\,$\beta$ according to Lemma \ref{lem:generator}. Then, we choose $\kappa_{i}$ according to Lemma \ref{lem:zero-dynamics}. Finally, we  choose $\rho_i$, $\tau_i$ to satisfy \eqref{eq:control-design-thm1}. 
	\end{remj}

In some case,  set $\W$  or at least its boundary might be known to us. Of course, we can still use the controller \eqref{ctr:adaptive} to solve this problem. But we can further construct a simpler controller based on the information of $\W$.  To this end, it is reasonable to introduce a new assumption to replace Assumption \ref{ass:zero-dynamics}.
\begin{assj}\label{ass:zero-dynamics-known}
	For each $i\in \mathcal{N}$, there exists a continuously differentiable function $W_{i\bar z}(\bar z_i)$ such that, for all $r_i\in \R$ and $w\in \W$, along the trajectory of system \eqref{sys:agent-error}, 
	\begin{align}\label{eq:zero-dynamics-known}
	\begin{split}
	&\underline{\alpha}_{i}(||{\bar z}_i||)\leq W_{i\bar z}({\bar z}_i)\leq {\bar \alpha}_i(||{\bar z}_i||) \\ 
	&\dot{W}_{i\bar z}\leq -\alpha_i(||\bar z_i||)+  \gamma_{ie}(e_i) e_i^2+ \gamma_{i\mu}(r_i)\mu_i^2
	\end{split}
	\end{align}
	for some known smooth functions $\underline{\alpha}_i$, ${\bar \alpha}_i$, $\alpha_i \in \mathcal{K}_\infty$, $\gamma_{ie}$, $\gamma_{i\mu}>1$ with $\alpha_i$ satisfying $\limsup_{s \to 0+}\frac{\alpha_i^{-1}(s^2)}{s}< +\infty$.
\end{assj}

In this case, we propose a reduced-order  controller: 
\begin{align}\label{ctr:adaptive-none}
u_i=\;&-\rho_i(\zeta_i,\,r_i) \zeta_i+\kappa_i(r_i)\eta_i\nonumber \\
\dot{\eta}_i=\;&-\kappa_i(r_i)\eta_i+u_i \\
\dot{r}_i=\;&-\alpha \nabla f_i(r_i)-\beta \sum\nolimits_{j=1}^{N}a_{ij}(r_i-r_j)-\sum\nolimits_{j=1}^{N}a_{ij}(v_i-v_j)\nonumber\\
\dot{v}_i=\;&\alpha \beta  \sum\nolimits_{j=1}^{N}a_{ij}(r_i-r_j) \nonumber
\end{align}
The optimal output consensus problem can be solved by this new controller as shown in the following theorem.
\begin{thmj}\label{thm:main-adaptive-non}
	Under Assumptions \ref{ass:graph}--\ref{ass:re-solution} and \ref{ass:zero-dynamics-known}, there exist positive constants $\alpha,\,\beta$ and positive smooth functions $\kappa_{i}(r_i)$, $\rho_i(\zeta_i,\,r_i)$ such that Problem \ref{prop} for multi-agent system \eqref{sys:agent} is solved by a distributed controller of the form \eqref{ctr:adaptive-none}. 
\end{thmj}
\pb  The proof is similar as that of Theorem \ref{thm:main}, and we only provide some brief arguments.  
 
First, by similar arguments as that in the proof of Lemma \ref{lem:zero-dynamics}, we can show that, for each $i\in \mathcal{N}$, there exist a smooth function $\kappa_i(r_i)>0$ and a continuously differentiable function $W_i(\hat z_i)$ such that, along the trajectory of system \eqref{sys:agent-error:translated-im},
	\begin{align}
	\begin{split}
	&\hat {\underline \alpha}_{i}(||\hat z_i||) \leq W_i(\hat z_i)\leq \hat {\bar \alpha}_i(||\hat z_i||) \\ 
	&\dot{W}_i\leq -||\hat z_i||^2 +\hat \gamma_{i\zeta}(\zeta_i,\,r_i)\zeta_i^2+\hat \gamma_{i\mu}(\mu_i,\,r_i)\mu_i^2
	\end{split}
	\end{align}
for some known smooth functions $\hat {\underline{\alpha}}_i,\,\hat {\bar \alpha}_i \in \mathcal{K}_\infty$,\,$\hat \gamma_{i\zeta},\,\hat \gamma_{i\mu}>1$. 

Next, we apply the changing supply functions technique to the $\hat z_i$-subsystem and conclude that, for any given smooth function $\hat \Delta_i(\hat z_i)>0$, there exists a continuously differentiable function $W_i^1(\hat  z_i)$ such that, along the trajectory of \eqref{sys:agent-error:translated-im},
\begin{align*}
&\hat {\underline \alpha}^1_{i}(||\hat   z_i||) \leq W_i^1(\hat z_i)\leq \hat {\bar \alpha}^1_i(||\hat  z_i||)\\
&\dot{W}_i^1\leq -\hat \Delta_i(\hat z_i)||\hat  z_i||^2+\hat \gamma^1_{i\zeta}(\zeta_i,r_i)||\zeta_i||^2+ \hat \gamma^1_{i\mu}(\mu_i,r_i)\mu_i^2
\end{align*}
for some known smooth functions $\hat {\underline\alpha}^1_i,\,\hat {\bar \alpha}^1_i \in \mathcal{K}_\infty$,\,$\hat \gamma^1_{i\zeta},\,\hat \gamma^1_{i\mu}>1$.

Let $\hat V_i(\hat z_i,\,\zeta_i)=W_i^1(\hat z_i)+\zeta_i^2$. By Lemma 11.1 in \cite{chen2015stabilization} and completing the square, one can obtain that
\begin{align*}
\dot{\hat V}_i
\leq\;& -[\hat \Delta_i(\hat z_i)-\check\phi_{i2}(\hat  z_i)]||\hat  z_i||^2+ [\hat \gamma^1_{i\mu}(\mu_i,\,r_i)+k^2_{i1}]\mu_i^2\\
& -[2b_i(w)\rho_i(\zeta_i,\,r_i)-\hat \gamma^1_{i\zeta}(\zeta_i,\,r_i)-\check\phi_{i1}(r_i)-\check\phi_{i3}(\zeta_i)-1]\zeta_i^2
\end{align*}
for some  known smooth functions $\check \phi_{i1}, \check \phi_{i2}, \check \phi_{i3}>1$. Letting $\hat \Delta_i(\hat z_i) \geq\check \phi_{i2}(\hat  z_i)+1$, $\rho_i(\zeta_i,\,r_i) \geq\frac{1}{2b_0}[\hat \gamma^1_{i\zeta}(\zeta_i,\,r_i)+\check \phi_{i1}(r_i)+\check \phi_{i3}(\zeta_i)+2]$ implies $\dot{\hat V}_i\leq -||\hat z_i||-\zeta_i^2+\hat c_{i\mu}\mu_i^2$ for some constant $\hat c_{i\mu}>0$. Then, the arguments of {\it Step 2} in the proof of Theorem \ref{thm:main} proceed as well and thus complete the proof.  
\pe

Since we are supposed to know the boundary of set $\W$, no adaptive component is needed in  controller \eqref{ctr:adaptive-none}. In this case, the rest parameters and functions can be derived in a similar way as mentioned in Remark \ref{rem:control}.

\begin{remj}\label{rem:embed}
	The controllers \eqref{ctr:adaptive} and \eqref{ctr:adaptive-none} are both composed of two parts constructed in two steps: optimal signal generator for problem \eqref{opt:main} and distributed partial stabilizer for the augmented system composed of \eqref{sys:composite-osg} and \eqref{sys:agent-error}. By this two-step procedure and dynamic compensator based feedback designs, the technical difficulties brought by nonlinearities, uncertainties and optimal requirements are successfully overcame.
\end{remj}

\begin{remj}\label{rem:ijrnc}
	Compared with relevant reference \cite{tang2018ijrnc}, the multi-agent system \eqref{sys:agent} is further subject to dynamic uncertainties. Moreover, the considered agents are nonlinearly parameterized with respect to uncertainties in contrast to the linear parameterized fashion in \cite{tang2018ijrnc}.  As the pure adaptive rules fail to solve this problem, a novel robust distributed controller has been developed to deal with the complicated uncertainties.
\end{remj}

\section{Simulation}\label{sec:simu}

In this section, we present two examples to illustrate the effectiveness of our designs.

\begin{figure}
	\centering
	\begin{tikzpicture}[shorten >=1pt, node distance=1.8cm, >=stealth',every state/.style ={circle, minimum width=0.3cm, minimum height=0.3cm}]
	\node[align=center,state](node1) {\scriptsize 1};
	\node[align=center,state](node2)[right of=node1]{\scriptsize 2};
	\node[align=center,state](node3)[right of=node2]{\scriptsize 3};
	\node[align=center,state](node4)[right of=node3]{\scriptsize 4};
	\path[->]  (node1) edge (node2)
	(node2) edge (node3)
	(node3) edge (node4)
	(node4) edge [bend right=22]  (node1);
	\path[<->] (node1) edge [bend right=20] (node3)
		(node2) edge [bend right=20]  (node4)
	;
	\end{tikzpicture}	
	\caption{Communication digraph $\mathcal G$ in our examples.}\label{fig:graph}
\end{figure}
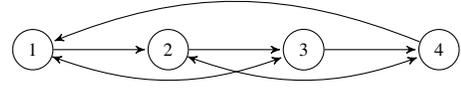

\begin{figure}
	\centering
	\includegraphics[width=0.44\textwidth]{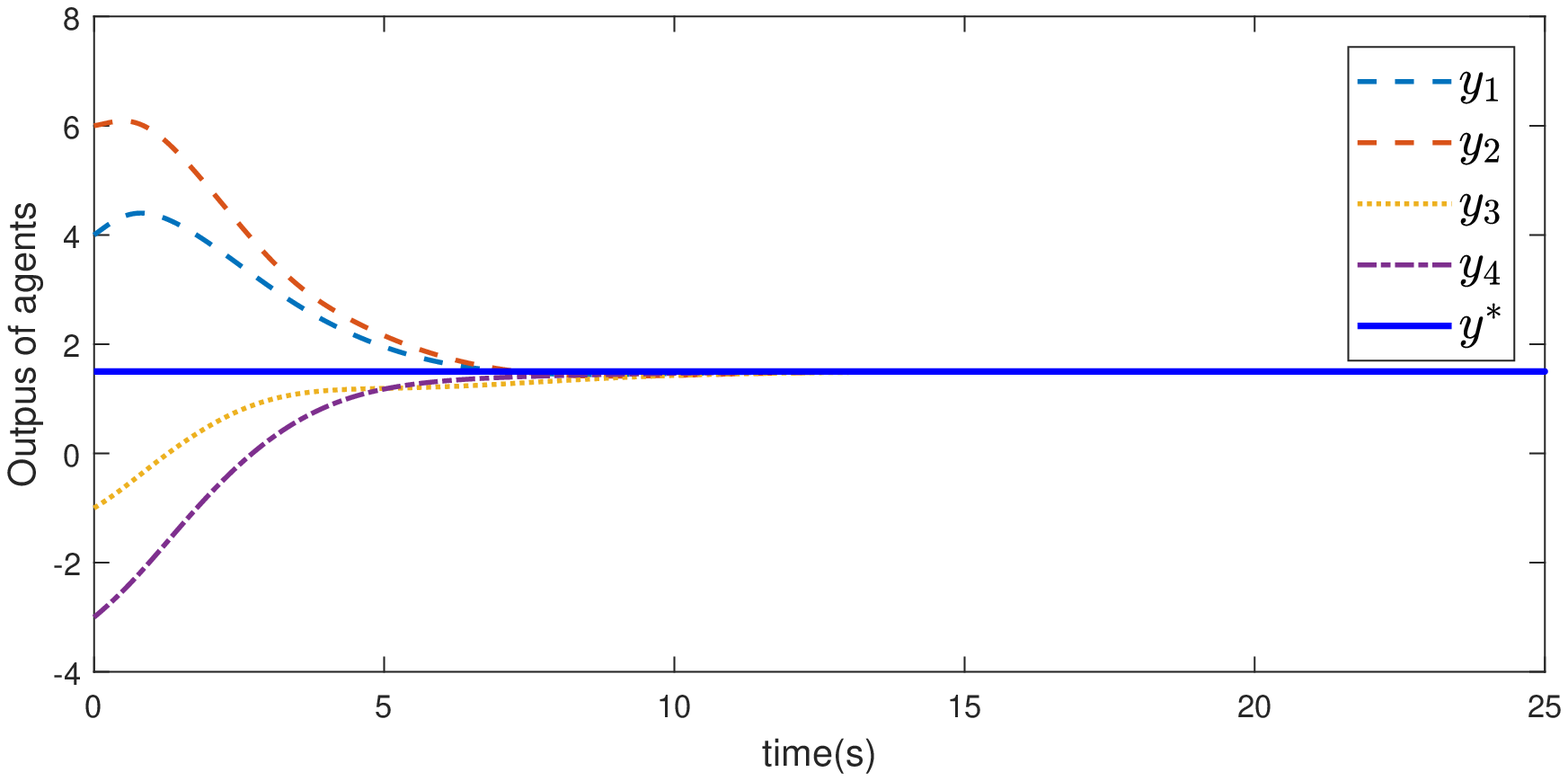}
	\caption{Profiles of agent outputs in Example 1. }\label{fig:simu1} 
\end{figure}

{\em Example 1} Consider a rendezvous problem \cite{ye2015model} for four single-link manipulators with flexible joints  as follows:
\begin{align}\label{sys:manipulator-exam}
\begin{split}
J_{i1} \ddot{q}_{i1}+M_igL_i \sin q_{i1}+ k_i(q_{i1}-q_{i2})&=0\\
J_{i2} \ddot{q}_{i2}-k_i(q_{i1}-q_{i2})&=u_i
\end{split}
\end{align}
where $q_{i1}, q_{i2}$ are the angular positions, $J_{i1},\, J_{i2}$ ar the moments of inertia, $M_i$ is  the total mass, $L_i$ is a distance, $k_i$ is a spring constant, and $u_i$ is the torque input. The communication digraph among these agents is depicted as Fig.~\ref{fig:graph} with unity edge weights with $\lambda_2=2$ and $\lambda_4=3$. 

To steer these manipulators to rendezvous at a common position that minimizes the aggregate distance from their starting position to this final position, we let $y_i=q_{i1}$ and take the cost functions as $f_i(y_i)= \frac{1}{2}||y_i-q_{i1}(0)||^2$ and $f(y)= \frac{1}{2}\sum_{i=1}^4 ||y-q_{i1}(0)||^2$ ($i=1,\,\dots,\, 4$). One can check that the optimal solution of the global cost function is $y^{\star}=\frac{1}{4}\sum_{i=1}^4 q_{i1}(0)$. To make this problem more interesting, we assume that $M_i=(1+w_{i1})M_{i0}$ and $L_i=(1+w_{i2})L_{i0}$ with nominal mass $M_{i0}$, nominal length $L_{i0}$, and uncertain parameters $w_{i1},\, w_{i2}$. 

Letting $x_i=\mbox{col}(q_{i1},\,\dot{q}_{i1},\,{q}^{(2)}_{i1},\,{q}^{(3)}_{i1})$, we rewrite system \eqref{sys:manipulator-exam} into the form \eqref{sys:agent} with  $w=\mbox{col}(w_{11},\,w_{12},\,\dots,\,w_{41},\,w_{42})$, $n_i=4$, $b_i(w)=\frac{k_i}{J_{i1}J_{i2}}$ and $g_i(x_i,\, w)=-x_{i3}[\frac{M_igL_i}{J_{i1}}\cos(x_{i1}) + \frac{k_i}{J_{i1}}+ \frac{k_i}{J_{i2}}]+ \frac{M_igL_i}{J_{i1}}({x}_{i2}^2-\frac{k_i}{J_{i2}})\sin(x_{i1})$. We can verify all assumptions in this paper and solve this problem according to Theorem \ref{thm:main}.

For simulations, we set $J_{i1}=1$, $J_{i2}=1$, $L_{i0}=1$, $M_i=1$, $k_i=1$ for simplicity and the uncertain parameters are randomly chosen such that $w_{i1},\, w_{i2}\geq 0$. Following the procedures in Lemma \ref{lem:generator} and Theorem \ref{thm:main}, we select $\alpha=1$, $\beta=15$ for the generator \eqref{sys:composite-osg} and $k_{i1}=1$, $k_{i2}=3$, $k_{i3}=3$, $\kappa_i(r_i)=1$, $\rho_i(\zeta_i,\,r_i)=\zeta_i^4+1$, $\tau_i(\zeta_i,r_i)=\rho_i(\zeta_i,\,r_i)\zeta_i^2$ for the controller \eqref{ctr:adaptive} with $1\leq i\leq 4$.  All initial conditions are randomly chosen and the simulation result is shown in Fig.~\ref{fig:simu1}, where the optimal rendezvous can be observed on $y^{\star}$. 

{\em Example 2} Consider another multi-agent system including two controlled FitzHugh-Nagumo dynamics \cite{murray2002mathematical}
\begin{align*}
\dot{z}_{i}=\;&-(1+w_{i3})cz_{i}+(1-w_{i4})bx_{i}\\
\dot{x}_{i}=\;&(1+w_{i6})x_{i}(a-x_{i})(x_{i}-1)-z_{i}+(1+w_{i5})u_i\\
y_{i}=\;&x_{i}, \quad  i=1,\,2  
\end{align*}
and two controlled Van der Pol oscillators \cite{khalil2002nonlinear}
\begin{align*}
\dot{x}_{i1}=\;&x_{i2}\\
\dot{x}_{i2}=\;&-(1+w_{i3})x_{i1}+(1+w_{i4})(1-x_{i1}^2)x_{i2}+(1+w_{i5})u_i\\
y_i=\;&x_{i1},\quad i=3,\,4
\end{align*}
with input $u_i$, output $y_i$, constants $a,\,b,\,c>0$, and unknown parameter $w_{ij}$.  Let $w=\mbox{col}(w_{13},\,w_{14},\,\dots,\,w_{44},\,w_{45})$. Clearly, all these agents are of the form \eqref{sys:agent}. 

We consider the optimal output consensus problem for this heterogeneous multi-agent system with more complicated cost functions  as ${f_1}(y) = (y-8)^2$, ${f_2}(y) = \frac{y^2}{80\ln {({y^2} + 2} )} + (y - 5)^2$,
${f_3}(y)=\frac{y^2}{{20\sqrt {y^2 + 1} }} + y^2$, ${f_4}(y) =  \ln \left( {{e^{ - 0.05{y}}} + {e^{0.05{y}}}} \right) + y^2$. Using the inequalities $0\leq\frac{1}{\ln(y^2+2)}\leq 1.5$, $0\leq\frac{1}{\sqrt{y^2+1}}\leq 1$, $-1\leq\frac{e^{0.05y}-e^{-0.05y}}{e^{0.05y}+e^{-0.05y}}\leq 1$, we can verify Assumption \ref{ass:convexity-strong} with $\underline{l}_i=1$ and $\bar{l}_i=3$ for $i=1,\,\dots,\, 4$. Furthermore, the global optimal point is $y^{\star}=3.24$ by numerically minimizing $\sum_{i=1}^4 f_i(y)$.

\begin{figure}
	\centering
		\includegraphics[width=0.44\textwidth]{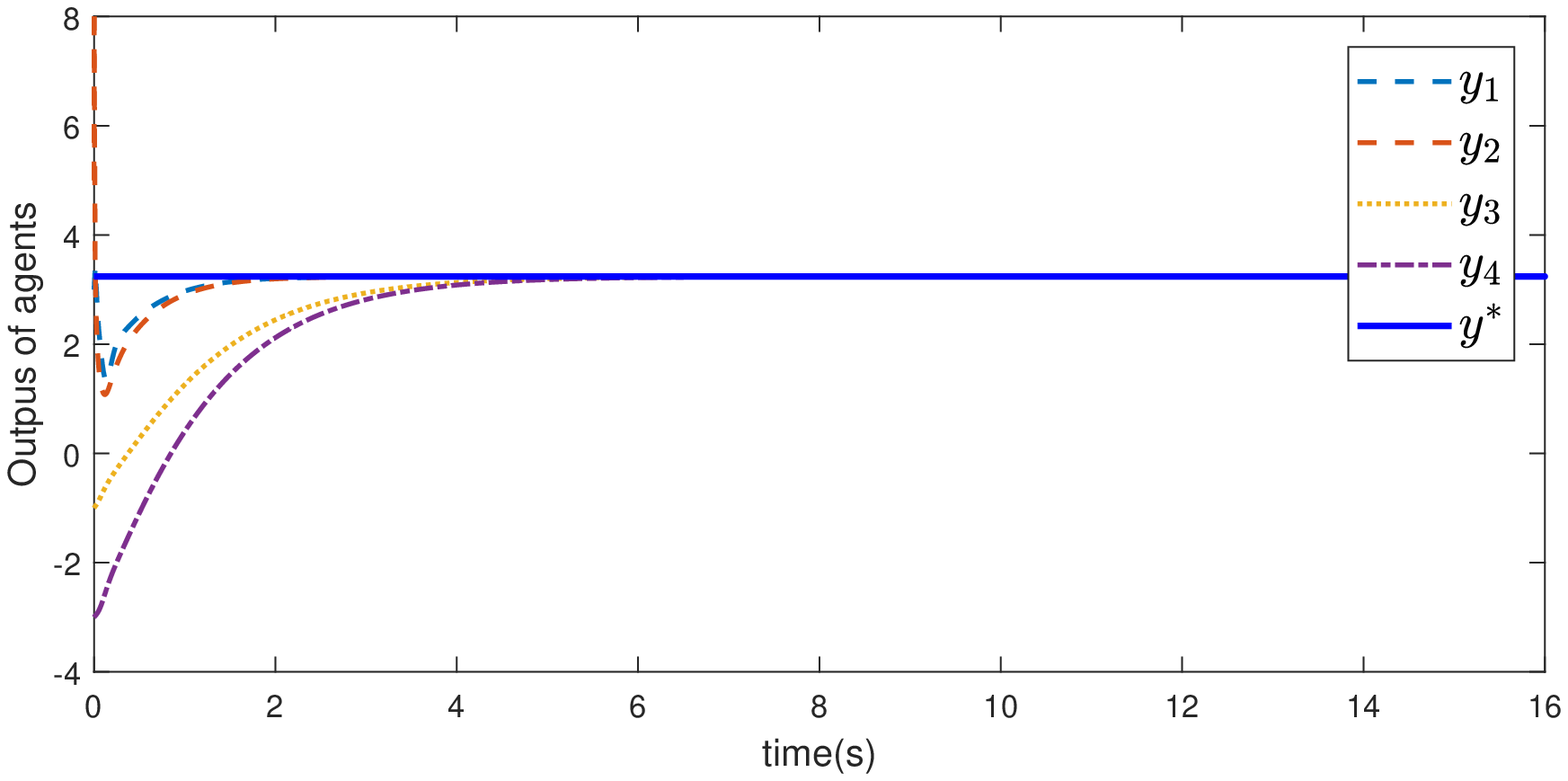} 
	\caption{Profiles of agent outputs in Example 2. }\label{fig:simu2} 
\end{figure}

Let $a=0.2$, $b=0.8$, $c=0.8$. The uncertain parameters are randomly chosen such that $w_{i3},\, w_{i5}\geq 0$ for $i=1,\,\dots,\,4$.  Without knowing the boundary of the compact set $\mathbb{W}$ containing these uncertainties, the controllers in \cite{tang2018acc} fail to solve the associated optimal output consensus problem. However, we can verify Assumptions \ref{ass:re-solution} and \ref{ass:zero-dynamics} for $i=1,\,2$ with ${z}^{\star}_i(s,\,w)=\frac{(1-w_{i2})b}{(1+w_{i1})c}s$, $W_{i\bar z}(s)=\alpha_{i}(s)=s^2$, $\gamma_{ie}(s)=\gamma_{i\mu}(s)=1$. Note that these two assumptions trivially hold for $i=3,\, 4$. According to Theorem \ref{thm:main}, the associated optimal output consensus problem can be solved by a distributed controller of the form \eqref{ctr:adaptive}. For simulations, we still use $\alpha=1$, $\beta=15$, and then choose $\rho_i(\zeta_i,\,r_i)=\zeta_i^4+r_i^4+1$, $\kappa_{i}(r_i)=r_i^4+1$, $\tau_i(\zeta_i,r_i)=\rho_i(\zeta_i,\,r_i)\zeta_i^2$ with $\zeta_i=x_{i}-r_i$ for $i=1,\,2$ and  $\rho_i(\zeta_i,\,r_i)=\zeta_i^4+r_i^4+1$, $\kappa_{i}(r_i)=r_i^4+1$, $\tau_i(\zeta_i,r_i)=\rho_i(\zeta_i,\,r_i)\zeta_i^2$ with $\zeta_i=x_{i1}-r_i+x_{i2}$ for $i=3,\,4$. All initial conditions are randomly chosen and the simulation result is shown in Fig.~\ref{fig:simu2}, where a satisfactory performance can be observed and the optimal output consensus is achieved on the optimal point $y^{\star}=3.24$.

\section{Conclusion}\label{sec:con}
We have studied an optimal output consensus problem  for a class of heterogeneous high-order nonlinear systems with both static and dynamic uncertainties. We proposed a two-step design scheme to convert it into two subproblems: optimal consensus for single-integrator multi-agent system and distributed partial stabilization of some augmented nonlinear systems. By adding a dynamic compensator to deal with the uncertainties, we constructed two distributed controls for this problem under standing conditions.  Our future works include the MIMO extension with time-varying digraphs.


\section*{Appendix.  Proof of Lemma \ref{lem:zero-dynamics}}

The proof is completed by successively using the changing supply functions technique \cite{sontag1995changing}.  

We first consider the case when $n_i\geq 2$. 
Under Assumption \ref{ass:zero-dynamics}, we apply the changing supply functions technique to  the $\bar z_i$-subsystem and conclude that, for any given $\bar \Delta_{i\bar z}(\bar z_i)>0$, there exists a continuously differentiable function $W_{i\bar z}^1(\bar z_i)$ satisfying   
$\bar {\underline \alpha}_{i\bar z}(||{\bar z_i}||) \leq W^1_{i\bar z}({\bar z_i})\leq \bar {\bar \alpha}_{i\bar z}(||{\bar z_i}||)$ and 
\begin{align*}
\dot{W}^1_{i\bar z} \leq -\bar\Delta_{i\bar z}({\bar z_i})||{\bar z_i}||^2+  \bar\sigma_{i\xi} \bar \gamma_{i\xi}^1({\xi_i})||{\xi_i}||^2+  \bar\sigma_{i\mu} \bar \gamma^1_{i\mu}(\mu_i,\,r_i)\mu_i^2
\end{align*}
for some known smooth functions $\bar{\underline\alpha}_{i\bar z}$,\,$\bar{\bar \alpha}_{i\bar z} \in \mathcal{K}_\infty$, $\bar \gamma^1_{i\zeta}$, $\bar \gamma^1_{i\mu}>1$ and unknown constants $\bar\sigma_{i\xi}, \bar \sigma_{i\mu}>1$. 

From the choice of $k_{ij}$, matrix $A_i^{\rm o}$ is Hurwitz. Then, there exists a unique positive definite matrix $\bar P_i$ satisfying  ${A_i^{\rm o}}^\top \bar P_i + \bar P_i A_i^{\rm o}=-3I_{m-1}$.  Let $W_{i\xi}^0(\xi_i)=\xi_i^\top \bar P_i \xi_{i}$. Its time derivative along the trajectory of \eqref{sys:agent-error} satisfies
\begin{align*}
\dot{W}_{i\xi}^0=\;& 2 \xi_i^\top \bar P_i[A_i^{\rm o} {\xi}_i+B_i^{\rm o} \zeta_i- E_i^{\rm o}  \mu_i]\\
\leq\;& - ||\xi_i||^2+||\bar P_iB_i^{\rm o}||^2 ||\zeta_i||^2+ ||\bar P_iE_i^{\rm o}||^2\mu_i^2
\end{align*}
By changing supply functions of $\xi_i$-subsystem, for any given $\bar \Delta_{i\xi}(\xi_i)>0$, there exists a continuously differentiable function $W_{i\xi}^1(\xi_i)$ 
satisfying 
$\bar {\underline \alpha}_{i\xi}(||{\xi_i}||) \leq W_{i\xi}^1({\xi_i})\leq \bar {\bar \alpha}_{i\xi}(||{\xi_i}||)$ and 
\begin{align*}
\dot{W}^1_{i\xi} \leq -\bar\Delta_{i\xi}({\xi_i})||{\xi_i}||^2+ \bar \gamma_{i\zeta}({\zeta_i})||{\zeta_i}||^2+  \bar \gamma_{i\mu}(\mu_i)\mu_i^2
\end{align*}
for some known smooth functions  $\bar {\underline \alpha}_{i\xi},\, \bar {\bar \alpha}_{i\xi} \in \mathcal{K}_\infty$, $\bar \gamma_{i\zeta}$,\,$\bar \gamma_{i\mu}>1$.

Let $W_{i\tilde z}(\tilde z_i)=W_{i\bar z}^1(\bar z_i)+\bar\sigma_{i\xi} W_{i\xi}^1(\xi_{i})$. Clearly, there exist functions $\tilde {\underline \alpha}_{i},\,\tilde {\bar \alpha}_i \in \mathcal{K}_\infty$ satisfying  $\tilde {\underline \alpha}_{i}(||\tilde  z_i||) \leq W_{i\tilde z}(\tilde  z_i)\leq \tilde {\bar \alpha}_i(||\tilde  z_i||)$.  Its time derivative along the trajectory of \eqref{sys:agent-error:translated-im} satisfies
\begin{align*}
\dot{W}_{i\tilde z}
\leq\;& -\bar\Delta_{i\bar z}({\bar z_i})||{\bar z_i}||^2-\bar\sigma_{i\xi}(\bar\Delta_{i\xi}({\xi_i})- \bar \gamma_{i\xi}^1({\xi_i}))||{\xi_i}||^2\\
&+\bar\sigma_{i\xi}\bar \gamma_{i\zeta}({\zeta_i})||{\zeta_i}||^2+  \bar\sigma_{i\xi}\bar \gamma_{i\mu}(\mu_i)\mu_i^2+\bar\sigma_{i\mu} \bar \gamma^1_{i\mu}(\mu_i,\,r_i)\mu_i^2
\end{align*}
Letting  $\bar\Delta_{i\bar z}({\bar z_i})>1$, $\bar\Delta_{i\xi}({\xi_i})> \bar \gamma_{i\xi}^1({\xi_i})+1$, $\tilde  \gamma_{i\mu}(\mu_i,\,r_i)> \bar \gamma_{i\mu}(\mu_i)+\bar \gamma^1_{i\mu}(\mu_i,r_i)$, and $\tilde \sigma_{i\mu}>\max\{\bar\sigma_{i\xi}, \, \bar\sigma_{i\mu} \}$ gives
\begin{align*}
\dot{W}_{i\tilde z}\leq -||\tilde z_i||^2+ \bar \sigma_{i\zeta} \bar \gamma_{i\zeta}(\zeta_i)||\zeta_i||^2+  \tilde \sigma_{i\mu} \tilde  \gamma_{i\mu}(\mu_i,r_i)\mu_i^2
\end{align*}
When $n_i=1$, the above property trivially holds for  $\tilde z_i= \bar z_i$.

Next, we apply the changing supply functions technique to $\tilde z_i$-subsystem and conclude that, for any given smooth function $\tilde \Delta_i(\tilde z_i)>0$, there exists a continuously differentiable function $W^1_{i\tilde z}(\tilde  z_i)$ satisfying $\tilde {\underline \alpha}^1_{i}(||\tilde  z_i||) \leq W_{i\tilde z}^1(\tilde  z_i)\leq \tilde {\bar \alpha}^1_i(||\tilde  z_i||)$ and
\begin{align*}
\dot{W}^1_{i\tilde z} \leq -\tilde \Delta_i(\tilde  z_i)||\tilde  z_i||^2+ \tilde \sigma_{i\zeta} \tilde  \gamma^1_{i\zeta}(\zeta_i)||\zeta_i||^2+  \tilde \sigma_{i\mu} \tilde  \gamma^1_{i\mu}(\mu_i,r_i)\mu_i^2
\end{align*}
for some known smooth functions $\tilde  {\underline\alpha}^1_i$,\,$\tilde  {\bar \alpha}^1_i \in \mathcal{K}_\infty$, $\tilde \gamma^1_{i\zeta}$,\,$\tilde  \gamma^1_{i\mu}>1$, and unknown constants $\tilde \sigma_{i\zeta}, \tilde \sigma_{i\mu}>1$.

Let $W_i(\hat z_i)= \tilde \ell_i   W_{i\tilde z}^1(\tilde  z_i)+\bar \eta_i^2  $ with $\tilde \ell_i>0$ to be specified later. Clearly, the first inequality in Lemma \ref{lem:zero-dynamics} holds. We take time derivative of $W_i$ along the trajectory of \eqref{sys:agent-error-proof} and have 
\begin{align*}
\dot{W}_i\leq\;& - \tilde \ell_i[\tilde \Delta_i(\tilde  z_i)||\tilde  z_i||^2- \tilde \sigma_{i\zeta} \tilde \gamma^1_{i\zeta}(\zeta_i)||\zeta_i||^2-\tilde \sigma_{i\mu} \tilde  \gamma^1_{i\mu}(\mu_i,r_i)\mu_i^2]\\
&+2 \bar \eta_i [-\kappa_i(r_i)\bar \eta_i+\tilde g_i(\bar z_i,\,\xi_i,\, \zeta_i,\, r_i,\, w)+\psi_i(r_i,\,w)\mu_i]
\end{align*}

Jointly with the inequalities \eqref{eq:growth11} and \eqref{eq:growth12}, we can bound the cross terms by completing the square and have
\begin{align*}
\dot{W}_i\leq\;& -[\tilde \ell_i\tilde \Delta_i(\tilde  z_i)-\frac{2\hat c_{ig} \hat \phi_{i2}(\tilde  z_i)}{b_0^2}]||\tilde  z_i||^2\\
&-[\kappa_i(r_i)-\frac{\hat \phi_{i1}(r_i)}{2}-\frac{1}{2}]\bar \eta_i^2\\
&+ [\tilde \ell_i\tilde \sigma_{i\zeta}\hat \gamma^1_{i\zeta}(\zeta_i)+ \frac{\kappa_i(r_i)}{b_0^2}+\frac{2\hat c_{ig} \hat \phi_{i3}(\zeta_i)}{b_0^2} ] ||\zeta_i||^2\\
&+ [2\hat \ell_{i\psi}\hat \phi_{i4}(r_i)+\tilde\ell_i\tilde \sigma_{i\mu}\hat \gamma^1_{i\mu}(\mu_i,\,r_i)]\mu_i^2
\end{align*}
Note that $\kappa_i(r_i)\geq \hat \phi_{i1}(r_i)+1$. Letting $\tilde \ell_i>\frac{2\hat c_{ig} }{b_0^2}+1$, $\tilde \Delta_i(\tilde  z_i)> \hat \phi_{i2}(\tilde  z_i)+1$,  $\hat \sigma_{i\zeta}>\tilde \ell_i\tilde \sigma_{i\zeta}+\frac{2\hat c_{ig}}{b_0^2}$, $\hat \sigma_{i\mu}>\tilde  \ell_i\tilde \sigma_{i\mu}+2\hat \ell_{i\psi}$, $\hat \gamma_{i\zeta}(\zeta_i,r_i)>\hat \gamma^1_{i\zeta}(\zeta_i)+ \kappa_i(r_i)+ \hat \phi_{i3}(\zeta_i)$, and $\hat \gamma_{ir}(\mu_i,r_i)>\hat \phi_{i4}(r_i)+ \hat \gamma^1_{i\mu}(\mu_i,\,r_i)$ 
implies the second inequality and thus completes the proof.

\bibliographystyle{IEEETran}
\bibliography{opt-high-nonlinear}

\end{document}